\documentclass[11pt]{amsart2000} 
\usepackage{amssymb}
\newtheorem{fact}{Fact}[section]
\newtheorem{ques}[fact]{Question}
\newtheorem{question}[fact]{Question}

\newtheorem{definition}[fact]{Definition}

\newtheorem{theorem}[fact]{Theorem}
\newtheorem{corollary}[fact]{Corollary}

\newtheorem{cor}[fact]{Corollary}
\newtheorem{prop}[fact]{Proposition}

\theoremstyle{definition}
\newtheorem{example}[fact]{Example}
\newtheorem{remark}[fact]{Remark}
\newcommand\T{{\mathbb T}}

\newcommand\Q{{\mathbb Q}}
\newcommand\Z{{\mathbb Z}}
\newcommand\ZZ{{\mathbb Z}}
\newcommand\R{{\mathbb R}}
\newcommand\N{{\mathbb N}}

\newcommand\J{{\mathbb J}}
\newcommand\co{{\mathfrak c}}
\newcommand\supp{\operatorname{supp}}

\begin{document}
\setlength{\unitlength}{0.01in}
\linethickness{0.01in}
\begin{center}
\begin{picture}(474,66)(0,0)
\multiput(0,66)(1,0){40}{\line(0,-1){24}}
\multiput(43,65)(1,-1){24}{\line(0,-1){40}}
\multiput(1,39)(1,-1){40}{\line(1,0){24}}
\multiput(70,2)(1,1){24}{\line(0,1){40}}
\multiput(72,0)(1,1){24}{\line(1,0){40}}
\multiput(97,66)(1,0){40}{\line(0,-1){40}}
\put(143,66){\makebox(0,0)[tl]{\footnotesize Proceedings of the Ninth Prague Topological Symposium}}
\put(143,50){\makebox(0,0)[tl]{\footnotesize Contributed papers from the symposium held in}}
\put(143,34){\makebox(0,0)[tl]{\footnotesize Prague, Czech Republic, August 19--25, 2001}}
\end{picture}
\end{center}
\vspace{0.25in}
\setcounter{page}{37}
\title{van Douwen's problems related to the Bohr topology}
\author{Dikran Dikranjan}
\thanks{The author was an invited speaker at the Ninth Prague Topological 
Symposium.}
\thanks{Work partially supported by the Research Grant of MIUR ``Nuove 
prospettive nella teoria degli anelli, dei moduli e dei gruppi abeliani''
2000}
\address{Dipartimento di Matematica e Informatica, Universit\`{a} di Udine\\ 
Via della Scienze 206, 33100 Udine, Italy}
\email{dikranja@dimi.uniud.it}
\thanks{Dikran Dikranjan,
{\em van Douwen's problems related to the Bohr topology},
Proceedings of the Ninth Prague Topological Symposium, (Prague, 2001),
pp.~37--50, Topology Atlas, Toronto, 2002}
\begin{abstract} 
We comment van Douwen's problems on the Bohr topology of the abelian
groups raised in his paper \cite{vD} as well as the steps in the solution
of some of them. 
New solutions to two of the resolved problems are also given. 
\end{abstract}
\maketitle

\section{Introduction}

Although fourteen years have passed since Eric van Douwen's untimely
death, his work is still of importance to general topologists. 
His questions (many of which remain open) continue to provide ideas and
inspiration to all of us. 
This survey collects recent results connected to some of the problems set
by van Douwen in his posthumous paper \cite{vD} on Bohr topology (see also
\cite[\S 4]{vDW}, \cite{vDP} for the list of 200 questions raised by van
Douwen in his publications with references to the relevant literature). 
The selection of the problems as well as the the information about their
present status, are due to author's personal taste and knowledge, without
any pretension for completeness.

The Bohr topology of an abelian group $G$ is {\em the largest totally
bounded group topology} on $G$. It can be realized as the initial topology
on $G$ with respect to the family of all homomorphisms of $G$ into the
circle group $\T$, i.e., the topology induced on $G$ under the canonical
embedding $G\hookrightarrow \T ^{Hom\;(G,\T)}$. In the sequel $G^{\#}$
will denote the group $G$ equipped with the Bohr topology. The completion
of $G^{\#}$ is known as {\em the Bohr compactification} of $G$ and denoted
by ${\rm b}G$ (defined otherwise with the property that every homomorphism
of $G$ into a compact group $K$ extends to a continuous homomorphism of
the group ${\rm b}G$ in $K$).

In \S 1 we recall some well-known properties of the Bohr topology and the
problems related to the topological properties (dimension, normality and
realcompactness) of $G^{\#}$. In \S 2 we outline the problems related to
the discrete sets in the Bohr topology. To the most relevant problem,
namely the homeomorphism problem (on whether the Bohr topologies of
abelian groups of the same cardinality must be homeomorphic), is
discussed in \S 4. 
The last section is dedicated to problems related to retracts of $G^{\#}$.

\subsection*{Notation and terminology} 

The symbols $\N$, $\Z$ and $\Q$ are used for the set of positive integers,
the group of integers and the group of rationals, respectively. The
circle group $\T$ is identified with the quotient group $\R/\Z$ of the
reals $\R$ and carries its usual compact topology. The cyclic group of
order $m$ is denoted by $\Z_m$ and $B_\kappa=\bigoplus_\kappa \ZZ_2$
denotes the Boolean group of size $\kappa$. The $p$-adic integers are
denoted by $\J_p$.

Let $G$ be a group. The set of torsion elements of $G$ is denoted by ${\rm
tor}(G)$ (it is a subgroup of $G$ when $G$ is abelian). For abelian groups
$G,H$ we denote by $Hom (G,H)$ the set of all homomorphisms $G\to H$ and
call a map $\ell:G\to H$ {\em linear} if $\ell(x)=h(x)+a$ for some $h\in
Hom (G,H)$, $a\in H$ and for every $x\in G$. Since
$Hom(B_\kappa,\T)=Hom(B_\kappa,\ZZ_2)$, a typical subbasic open
neighbourhood of $0$ in $B_\kappa^\#$ is $\ker \chi$, for a homomorphism
$\chi: B_\kappa \rightarrow \ZZ_2$.

The symbol $\co$ stands for the cardinality of the continuum, so
$\co=2^{\aleph_0}$. For undefined terms see \cite{E,Fu,HR,J}.

\section{Some properties of the Bohr topology}

We list below some well known properties of the Bohr topology. 
\begin{itemize}
\item 
Every linear map $\ell:G^\#\to H^\#$ is continuous (i.e., the Bohr
topology is {\em functorial}). 
\item 
For every subgroup $H$ of $G$:
\begin{itemize}
\item 
$H^\#$ is a {\em topological} subgroup of $G^\#$; 
\item 
$(G/H)^\#$ carries the {\em quotient} topology of $G^\#/H$ (i.e., the map
$G^\#\to (G/H)^\#$ is open);
\item
if $G=H\times L$, then $G^\#=H^\#\times L^\#$. 
\end{itemize}
\item 
If $G$ has exponent $m$ (that is, $(\forall x\in G)\; mx=0$), then $G^\#$
carries the profinite topology (i.e., the finite-index subgroups of $G$
form a local base at 0).
\item 
The weight and the character of $G^\#$ coincide with $2^{|G|}$
(so that an infinite $G^\#$ is never metrizable).
\end{itemize}

\subsection{Dimension, normality and compactness}

van Douwen proved that $G^\#$ is zero-dimensional for every abelian group
$G$ (\cite[Theorem 4.8]{vD}) and asked whether it is strongly
zero-dimensional \cite[Question 4.10]{vD}. A positive answer was obtained
by Shakhmatov \cite{S} in the following more general setting. He proved
that every zero-dimensional totally bounded group is strongly
zero-dimensional. A further generalization was obtained by Hern\'andez
\cite{Her}, who proved that for every locally compact abelian group $G$
the Lebesgue covering dimensions of the groups $G$ and $G\sp +$ coincide
(here $G^+$ denotes the group $G$ equipped with its Bohr topology, defined
as in the discrete case as the initial topology on $G$ with respect to the
family of all {\em continuous} homomorphisms of $G\to \T$). For discrete
$G$ this gives another solution of van Douwen's problem.

van Douwen proved that $G^\#$ is not paracompact in cases of interest: 

\begin{fact} {\rm \cite[Fact 4.11]{vD}}
$G^\#$ is paracompact iff $G^\#$ is collectionwise normal iff $G$ is
countable. 
\end{fact} 

He asked whether $G^\#$ is (always) normal (\cite[Question 4.10]{vD}).
Trigos-Arrieta \cite{Tr2} showed that one can add also ``normal'' in the
above chain of equivalences, i.e., $G^\#$ is normal iff $G$ is
countable. 
This result is remarkable, since it gives an easy ``uniform'' example of
non-normal group topology on every uncountable abelian group. Indeed, the
known examples of non-normal Hausdorff topological groups were the
uncountable powers of $\Z$ or free topological groups \cite{HR}, i.e.,
groups with rather stringent conditions on the algebraic structure.

Let us discuss now the degree of compactness of $G^\#$. On one hand it is
totally bounded, so close to being compact in this sense. On the other
hand, some of the typical properties of the compact groups are not shared
by $G^\#$. It is well known that the compact groups are dyadic spaces
(\cite{Iv,Kuz}), so they have plenty of non-trivial convergent
sequences. Flor \cite{Flor} proved that the only convergent sequences in
$G^\#$ are the trivial ones (so $G^\#$ is never sequential). This follows
also from the following more general fact, proved by Glicksberg \cite{Gl}:
the only compact sets in $G^\#$ are finite. Flor \cite{Flor} proved also
that actually no non-trivial sequence in $G^\#$ converges to a point of
${\rm b}G$. van Douwen noticed that this follows from the following
stronger property of $G^\#$: {\em if $A\subseteq G^\#$ is infinite and
$p\in {\rm b}G$, then there exists a neighbourhood $U$ of $p$ in 
${\rm b}G$ such that $|A\setminus U|=|A|$} \cite[Theorem 1.1.3(a)]{vD}.

As far as other weaker forms of compactness are concerned, it was proved
by Comfort and Saks \cite{CSa} that $G^\#$ is never pseudocompact when $G$
is infinite. van Douwen strengthens this fact by showing that $G^\#$ is
never a Baire space, when $G$ is infinite \cite[Corollary 4.7(b)]{vD}.

A. Dow (cf.\ \cite[Fact 4.15]{vD}) noticed that $G^\#$ is not realcompact
when $|G|$ is Ulam-measurable. 
This motivated 

\begin{question} 
{\rm \cite[Question 4.16]{vD}} 
Is $G^\#$ realcompact if $|G|$ is not Ulam-measurable? 
\end{question}

This question was answered positively by Comfort, Hern\'andez and
Trigos-Arrieta (\cite[Theorem 3.3]{CHT1}). Moreover, they extended this
result to the locally compact case by proving that for a locally compact
abelian group $G$ the group $G^+$ is realcompact if and only if $G$ is
realcompact \cite[Theorem 3.8]{CHT1}.

\section{The relatively discrete sets}

van Douwen proved that $G^\#$ contains relatively discrete closed sets of
size $\kappa$ for every cardinal $\kappa\leq |G|$ that is not 
Ulam-measurable (cf.\ \cite[p. 1073]{vD}, this follows from Theorem
\ref{vD_Thm}). This is why he asked whether this can be done for {\em all}
cardinals $\kappa\leq |G|$, and in particular:

\begin{question}\label{QuesD} {\rm \cite[Question 4.14]{vD}, \cite[\S4,
Question 83]{vDW}}
Does $G^\#$ contain relatively discrete closed sets of size $|G|$? 
\end{question} 

The following strong positive answer was obtained by Hart and van Mill
\cite{HvM}: 

\begin{theorem}
Every subset $A$ of $G^\#$ has a relatively discrete subset $D$ such that 
$|D|=|A|$ and $D$ is closed in $G^\#$.
\end{theorem} 

The following question of van Douwen is obviously motivated by Question
\ref{QuesD} (note that the groups $G^\#$ have plenty of relatively
discrete sets, cf.\ \ref{vD_Thm}).

\begin{question}\label{Ques3} {\rm \cite[Question 4.17]{vD}, \cite[\S
4.3, p. 20, Question 85]{vDW}}
If $G^\#$ is infinite, does $G^\#$ have a relatively discrete subset that
is not closed ? 
\end{question} 

It was proved by Hart and van Mill \cite[Example 0.2]{HvM} that every
Boolean group $B_\kappa$ admits a relatively discrete subset that is not
closed. Namely, the set $[\kappa]^2$ has this property. In other words,
the subset $\mathcal{D}_\kappa=\{0\}\cup [\kappa]^2$ of $B_\kappa^\#$
carrying the induced topology has a relatively discrete non-closed subset
(namely, $[\kappa]^2$). 

In 1991, Ursul \cite{U} answered positively Question \ref{Ques3}. 
Later Protasov \cite{P} gave a short and clear proof of a more general 
result: If $A$ is an infinite totally bounded subset of a topological
group $G$, then the set $A·A^{-1}$ contains a discrete non-closed subset.

The following theorem of K. Kunen and W. Rudin \cite{KR} provides, in
particular, another positive answer to Question \ref{Ques3} (just take the
set $(A-A)\setminus \{0\}$).

\begin{theorem} 
Every $G^\#$ contains a relatively discrete set $A$ of size $|G|$ such that:
\begin{itemize}
\item[(a)] $A$ is $C^*$-embedded in ${\rm b}G$, 
\item[(b)] 0 is the only limit point of $A-A$,
\item[(c)] $A+A$ has no limit points in $G^\#$ if $(G:G[2])=|G|$.
\end{itemize}
\end{theorem}

An alternative approach to Question \ref{Ques3} is possible due to the
above mentioned example of Hart and van Mill \cite{HvM}. 
Clearly, the groups $G^\#$ that contain a copy of $\mathcal{D}_\omega$ are
settled by that example. 
Therefore, it remains to apply the following result proved in \cite{D}: 

\begin{prop}\label{NEW} 
Let $G$ be an abelian group, let $\kappa$ be an infinite cardinal and let
$f:\kappa\to G$ be an arbitrary function. 
Then the map $\mu_f:\mathcal{D}_\kappa \to G^\#$ defined by $\mu_f(0)=0$
and $\mu_f(\alpha, \beta)=f(\alpha)-f(\beta)$ is continuous. 
Moreover, if the family $\{f(\alpha): \alpha <\kappa\}$ is independent,
then $\mu_f$ is an embedding. 
\end{prop}

The only cases left out by the embedding part of Proposition \ref{NEW}
are the groups $G$ that do not contain an infinite direct sum of
non-trivial subgroups. 
In such a case $G$ either contains a copy of $\Z$, or a copy of the
Pr\"ufer group $\Z(p^\infty)$. 
Hence one can argue with $G=\Z$ and $G=\Z(p^\infty)$. In both cases one
can easily construct the desired set $A$ as follows. 
Let $B=\{b_n\}$ be a countable subset of an abelian group $G$ such that: 
\begin{itemize}
\item[(i)] 
$0\not\in B$ and $B\cap (B+B)=\emptyset$; 
\item[(ii)] 
there exists a metrizable totally bounded group topology $\tau$ on $G$
such that $b_n\to 0$ in $\tau$.
\end{itemize}
Then $A=(B-B)\setminus \{0\}$ is a relatively discrete non closed set of
$G^\#$.
Indeed, it is easy to see that $0$ belongs to the closure of $A$ in any
totally bounded topology of $G$. 
On the other hand, $A$ is relatively discrete in $(G,\tau)$, hence in
$G^\#$ too. 
In the case of $G=\Z$ (ii) can be ensured if $b_{n+1}/b_n$ tends to
infinity, according to \cite{ZP}.
To ensure also (i) take $b_n=n!$. For $G=\Z(p^\infty)$ the sequence
$b_n=p^{-n!}+\Z$ works.

One can prove that the map $\mu_f$ in Proposition \ref{NEW} is generic in
the following sense: 
if $\kappa>\co$ and $h:\mathcal{D}_\kappa\to H^\#$ is continuous with 
$h(0)=0$, then there exists an infinite $S\subseteq \kappa$ and 
$f:S\to D(H)$ such that $h$ coincides with $\mu_f$ when restricted to
$\{0\}\cup [S]^2$ (\cite{D}). 

\subsection{Relatively discrete, $C$-embedded sets}

van Douwen \cite[Theorem 1.1.3 (a)]{vD} showed that 

\begin{theorem}\label{vD_Thm} 
Every subset $A$ of $G^\#$ has a relatively discrete subset $D$ such that
$|D|=|A|, D$ is $C$-embedded in $G^\#$ and $C^ *$-embedded in ${\rm b}G$. 
\end{theorem}

An alternative proof of this result was given by Galindo and S. Hern\'
andez \cite{GH1}. The following extension of this result to the locally
compact case was obtained later by the same authors \cite{GH}. Let $G$ be
an abelian $\mathcal{L}_\infty$-group in the sense of \cite{V} (i.e., the
topology of $G$ is the intersection of a non-increasing sequence of
locally compact group topologies). For a subset $A$ of a group $G$ they
denote by $B(A)$ the minimum size of a family of bounded sets covering $A$
(a subset $B$ of a topological group $G$ is bounded in this sense if for
every neighbourhood $U$ of $0 \in G$ there exists a finite subset $F$ of
$B$ and a natural number $n$ such that $B$ is contained in
$F+\underbrace{U+\ldots+U}\sb {n}$ \cite{He}). The authors prove that
every subset $A$ of $G$ has a relatively discrete subset $B$ with
$|B|=B(A)$ and such that $B$ is $C$-embedded in $G\sp +$ and
$C^*$-embedded in ${\rm b}G$ (this is the completion of $G^+$ as in the
discrete case). When $G$ is discrete, this obviously gives Theorem
\ref{vD_Thm}.

\section{The homeomorphism problem}

The following question of van Douwen \cite[Question 515]{C1} is motivated
by the fact that none of the properties described in \S 2 depend on the
{\em algebraic} structure of the group $G$.

\begin{question}\label{hom_q} 
If $G$ and $H$ are discrete abelian groups of the same cardinality are
then $G^\#$ and $H^\#$ homeomorphic as topological spaces?
\end{question}

The question cannot be found explicitly in \cite{vD}. According to the
reference of \cite{CHT3}, it was posed in van Douwen's letters to W.
Comfort of 30 June 1986 and 9 May 1987.

\subsection{The first non-homeomorphism theorems}

A negative answer to Question \ref{hom_q} was given by Kunen \cite{K} in 1996: 

\begin{theorem}\label{K1} {\rm \cite{K}} 
$(\bigoplus_\omega \ZZ_p)^\#$ and $(\bigoplus_\omega \ZZ_q)^\#$ are not
homeomorphic for primes $p\ne q$.
\end{theorem}

The proof of this theorem makes use of a Ramsey-style partition property
of sequences $X=\langle x_s: s\in [\omega]^n\rangle$ in $\bigoplus_\omega
\ZZ_q$ (namely, every sequence splits into a sum of finitely many linearly
independent {\em normal forms}, cf.\ \cite{K}). This property is
established by means of Ramsey ultrafilters. Although Ramsey ultrafilters
exist under CH \cite{J}, the proof does not depend on CH.

Another solution was obtained about the same time by Watson and the
author \cite{DW}: 

\begin{theorem}\label{Th_DW} 
$B_\kappa^\#$ and $(\bigoplus_\kappa \ZZ_3)^\#$ are not homeomorphic for
$\kappa>2^{2^\co}$. 
\end{theorem}

A different partition theorem was used in \cite{DW}, based on splitting of 
the {\em supports} of the members of $X=\langle x_s: 
s\in [\kappa]^n\rangle$ in a many-variable Delta system lemma style
(essentially contained in \cite{W}).
Two features should be mentioned in the comparison of the partition
theorems used in \cite{K,DW}.
First, the ``elementary forms'' appearing in these partition theorems are
different, secondly, the combinatorics behind the proof is different too. 
Indeed, while the proof presented in \cite{K} works for countable groups
and fully exploits the fact that the {\em codomain} is a vector space
over a finite field, the proof in \cite{DW} requires larger groups, but
allows for a substantial generalization (see \S \ref{SecST}). 

It should be mentioned that this approach gives more than just
non-homeomorphism theorems. 
More precisely, the following theorem was proved by Kunen.

\begin{theorem}\label{K_THM} {\rm \cite{K}} 
If $f:\{0\}\cup [\omega]^4\to (\bigoplus_\omega\Z_q)^\#$ is a continuous
map, where $\{0\}\cup [\omega]^4$ carries the topology induced by
$(\bigoplus_\omega \Z_p)^\#$, then $f$ is constant on $[S]^4$ for
some infinite $S\subseteq \omega$.
\end{theorem}

Clearly, Theorem \ref{K1} follows immediately from Theorem \ref{K_THM}. 

\begin{remark}\label{REM}
\mbox{}
\begin{itemize}
\item[(a)] 
An analogous result holds for $\kappa>2^{2^{\co}}$ and continuous maps 
$f:\{0\}\cup [\kappa]^4\to (\bigoplus_\kappa\Z_3)^\#$, where $\{0\}\cup
[\kappa]^4$ carries the topology induced by $B_\kappa^\#$ \cite{DW}. 
Here $\kappa>\omega$ is needed only when 
$|\supp f(\alpha, \beta, \gamma, \delta)|$ is unbounded: 
{\em if $f:\{0\}\cup [\omega]^4\to (\bigoplus_\omega\Z_q)^\#$ is a
continuous map, where $|\supp f(\alpha, \beta, \gamma, \delta)|$ is
bounded, then $f$ is constant on $[S]^4$ for some infinite $S\subseteq 
\omega$}
\cite{DW}. 
Consequently, every continuous map $B_\kappa^\#\to (\bigoplus_\omega
\ZZ_3)^\#$ is constant on some infinite subset of $B_\kappa$ hence cannot
be a homeomorphism.
\item [(b)] 
It was shown later in \cite{D} that the group $\;\bigoplus_\kappa \ZZ_3$ 
in Theorem \ref{Th_DW} can be replaced by {\em any abelian group without
infinite Boolean subgroups } (cf.\ Corollary \ref{corXXX}). 
\end{itemize}
\end{remark} 

According to Proposition \ref{NEW}, there is a continuous one-to-one map 
$\mathcal{D}_\kappa\to G^\#$, with $\kappa=|G|$ for every infinite abelian
group $G$. 
This shows that one cannot argue with doubletons in Theorem \ref{K_THM}
or in Remark \ref{REM} (for $\kappa=\omega$ and $G=\bigoplus_\omega \Z_3$
this example was given in \cite{K}).

\subsection{Continuous maps in the Bohr topology}\label{SecST}

Theorem \ref{K_THM} and Remark \ref{REM} suggest that continuous maps
$f:B_\kappa^\#\to (\bigoplus_\kappa\Z_3)^\#$ are ``frequently constant''.
It is natural to relate this phenomenon to the fact that
$Hom(L,\bigoplus_\omega\Z_3)=0$ for every subgroup $L$ of $B_\kappa$ (so
every {\em linear} map $L\to \bigoplus_\omega\Z_3$ is constant). Hence a
general approach to Question \ref{hom_q} can be based on the following
idea: prove that {\em continuous} maps $f:G^\#\to H^\#$ are always {\em
linear} on some non-singleton (even infinite) subset of $G$. If
$Hom(L,H)=0$ for every subgroup $L$ of $G$ (this entails $G$ is torsion),
then one obtains a non-homeomorphism theorem in the framework of Question
\ref{hom_q}.
Let us give the relevant definition in a more precise form: 

\begin{definition} 
For abelian groups $G,H$ and $A\subseteq G$ a map $f: A\to H$ {\em is
linear on $A$} if $f$ coincides on $A$ with the restriction of a linear
map $\ell:\mbox{span}(A)\to H$. 
\end{definition} 

The next example, due to Comfort, Hern\'andez, Trigos-Arrieta \cite{CHT3}
(see also \ref{SecCHT}), shows that such an ``approximation'' by linear
maps is not available in general:

\begin{example} 
$\Q^\#$ is homeomorphic to 
$((\Q/\Z)\times \Z)^\#=(\Q/\Z)^\#\times \Z^\#$, so there exists an
embedding $j: (\Q/\Z)^\#\hookrightarrow \Q^\#$. 
Now $Hom(A, \Q)=0$ for every subgroup $A$ of $\Q/\Z$, so every linear map
$\Q/\Z\to \Q$ is constant. 
Thus, $j$ is never linear on any non-singleton $S\subseteq \Q/\Z$. 
\end{example}

\subsubsection{The Straightening Law for Boolean subgroups} 

The partition property of sequences $X=\langle x_s: s\in
[\kappa]^4\rangle$ in $\bigoplus_\omega \ZZ_3$ used in \cite{DW} is
modified in \cite{D} replacing the group $\bigoplus_\kappa \Z_3$ by a
direct sum of copies of an arbitrary countable abelian group. This
partition theorem is applied to obtain the following theorem in the
framework of the above mentioned programme of ``approximation of
continuous maps by linear ones''. For the sake of brevity we call it {\em
Straightening Law} in the sequel.

\begin{theorem}\label{ST} 
Let $H$ be an arbitrary abelian group. 
If $\kappa>2^{2^{\co}}$ and $f:\{0\}\cup [\kappa]^4\to H^\#$ is a
continuous map, where $\{0\}\cup [\kappa]^4$ carries the topology induced
by $B_\kappa^\#$, then $f$ is linear on $\{0\}\cup [S]^4$ for some
infinite $S\subseteq \kappa$. 
\end{theorem}

Consequently, either $f$ is constant on $[S]^4$, or $f$ coincides on
$[S]^4$ with the restriction of an {\it injective} linear map
$\ell:span(S)\to H$ (for some {\em smaller} $S$).

\begin{corollary} 
Let $H$ be an arbitrary abelian group. If $\kappa>2^{2^{\co}}$ and
$f:B_\kappa^\#\to H^\#$ is a continuous map, then there exists an
infinite $S\subseteq \kappa$ such that $f$ is linear on 
$S+S=\{0\}\cup [S]^2$.
\end{corollary}

\begin{cor}\label{corXXX} 
There exist no continuous finite-to-one maps $f:G^\#\to H^\#$ if the
abelian group $G$ contains Boolean subgroups of size $>2^{2^{\co}}$
and the abelian group $H$ admits no infinite Boolean subgroups. 
\end{cor}

This leaves open the following

\begin{ques}\label{ques1} 
For which infinite cardinals $\kappa$ does there exist a continuous 1-1
map from $B_\kappa^\#$ to any torsion-free group $H^\#$? 
\end{ques}

Such a $\kappa$ must be $\leq 2^{2^{\co}}$ according to Corollary
\ref{corXXX}. Let us note that also that it suffices the take
$H=\bigoplus_\omega \Q$.

It follows from the Corollary \ref{corXXX} that the Bohr topology can {\em
detect} Boolean subgroups in the following sense: {\em If $G$ and $H$ are
abelian groups such that the Bohr topologies of the (discrete) powers
$G^{2^{2^{\co}}}$ and $H^{2^{2^{\co}}}$ are homeomorphic, then
$G$
contains non-trivial Boolean subgroups iff $H$ does}. Note that $|G|=|H|$
need not hold and the groups may be finite.

\subsubsection{The Straightening Law in the general case} 

It is shown in \cite{Dp} that the Straightening Law holds (with
appropriate modifications) also in the more general context when the
Boolean group $B_\kappa$ is replaced by a vector space over the finite
field $\Z/p\Z$. Therefore, the Bohr topology can measure, roughly
speaking, also $p$-subgroups of the abelian groups. Let us note that
Theorem \ref{K_THM} can be considered as a first instance of the
Straightening Law in this context, since in that situation the constant
maps are precisely the restrictions of the linear ones.

The proofs in \cite{Dp} are based on an appropriate combination and
development of ideas from \cite{D} and \cite{K} that lead to a partition
theorem (in the spirit of \cite{K}, improving the ``independence'' property
of the normal forms) for `sequences' $X=\langle x_s: s\in
[\kappa]^n\rangle$ in direct sums $\bigoplus_\kappa K$, where $K$ is an
arbitrary countable abelian group. The price to pay is to increase the
cardinality asking $\kappa>\beth_{n-1}$. Indeed, easy examples show that
countable sequences and spaces cannot suffice in the case of vector spaces
over a countably infinite field. Let us note also that the partition
theorem from \cite{D}, easily extendible to $n$-ary sequences $X=\langle
x_s: s\in [\kappa]^n\rangle$ in direct sums $\bigoplus_\kappa K$, cannot
help either. Indeed, the normal forms in these partition theorems need not
be independent as in \cite{K}, so such partition theorem cannot be used to
obtain a proof of the Straightening Law following the argument in 
\cite{D}. 
It is shown in \cite{Dp} that with some additional care such a kind of
independence can be achieved.

From now on $p$ will be an arbitrarily fixed prime number and $V_\kappa$
will be the vector space of size $\kappa$ over $\Z_p$. In these terms the
Straightening Law in the general case can be announced as follows: {\em
For $\kappa > \beth_{2p-1}$ and a continuous map $\pi:\{0\}\cup
[\kappa]^{2p}\to H^\#$, where $\{0\}\cup [\kappa]^{2p}$ is equipped with
the induced from $V_\kappa^\#$ topology, there exists an infinite
$S\subseteq \kappa$ such that $\pi$ is linear on $[S]^{2p}$} \cite{Dp}. In
analogy to the case of Boolean subgroups, one can replace $[S]^{2p}$ by
$[S]^{p}$ and assume that the linear map is either constant or injective.

\begin{theorem}\label{lastT} 
If $\kappa> \beth_{2p-1}$ and $f:V_\kappa^\#\to H^\#$ is a continuous map,
then there exists an infinite $S\subseteq \kappa$ such that $f$ is linear
on $[S]^p$. 
\end{theorem}

This theorem implies that for every continuous self-map
$f:(\bigoplus_\kappa \Z_{p^2})^\#\to (\bigoplus_\kappa \Z_{p^2})^\#$ with
$\kappa> \beth_{2p-1}$ there exists a subset $A$ of $\bigoplus_\kappa
\Z_{p^2}$ such that $pA$ is infinite and $f$ is linear on $[pA]^p$.

Theorem \ref{lastT} has also the following immediate consequence. 

\begin{corollary}\label{MainTT} 
If the abelian group $H$ contains no copies of $V_\omega$, then there
exists no continuous finite-to-one map
$\pi:(\bigoplus_{\kappa}\Z_p)^\#\to H^\#$ for $\kappa>\beth_{2p-1}$. 
\end{corollary}

Now we see that there exists no continuous finitely many-to-one map $G^\#
\to H^\#$, when $H$ is torsion-free and $G$ is an abelian group with
large $p$-subgroups for some prime number $p$. 

\begin{cor}\label{CorXX} 
Let $p$ be a prime number. 
There exist no finitely many-to-one continuous maps $f:G^\#\to H^\#$ if
$G$ contains $p$-subgroups of size $>\beth_{2p-1}$ and $H$ contains no
copies of $V_\omega$. 
\end{cor}

Following \cite{TY} we say that an abelian group $H$ is {\em almost
torsion-free} if $H$ contains no copies of $V_\omega$ for every prime $p$.
Corollary \ref{CorXX} implies that if $G^\#$ and $H^\#$ are homeomorphic
and $H$ is almost torsion-free, then $|{\rm tor}(G)|\leq \beth_\omega$.

These results give many pairs of groups $G$ and $H$ non-homeomorphic in
the Bohr topology. For example, when $H$ is almost torsion-free with
$|H|=\kappa>\beth_{2m-1}$, $L$ is arbitrary with $|L|\leq\kappa$, and
$G=L\times(\bigoplus_\kappa \Z_m)$.

\begin{corollary} 
Let $p$ be a prime number and let $G$, $H$ be abelian groups such that the
powers $H^\kappa$ and $G^\kappa$ are homeomorphic in the Bohr topology for
some infinite cardinal $\kappa\geq \beth_{2p-1}$.
Then $G$ contains non-trivial $p$-subgroups iff $H$ does. 
\end{corollary}

\begin{corollary} 
If $\kappa\geq \beth_\omega$ and $H^\kappa$ and $G^\kappa$ are
homeomorphic in the Bohr topology, then for every prime $p$, $G$ contains
non-trivial $p$-subgroups iff $H$ does. 
In particular, $G$ is torsion-free iff $H$ is torsion-free. 
\end{corollary}

\subsubsection{Stronger form of the Straightening Law}

Let us discuss now the limits of the Straightening Law. We do not know
whether the following stronger forms of the Straightening Law hold true.
The first one is determined by {\em smaller domain}:

\newtheorem*{SSL}{Strong Straightening Law}
\begin{SSL}
For every prime number $p$ and for every continuous map
$\pi:V_\omega^\#\to H^\#$ there exists an infinite set
$Z\subseteq\omega$ such that $\pi$ is linear on $[Z]^p$.
\end{SSL}

It obviously implies that for no almost torsion-free abelian group $H$
there is a 1-1 map $(\Q/\Z)^{(\omega)}\to H$ continuous in the Bohr
topology. Moreover, under this conjecture, a countable $p$-group $G^\#$
embeds into some torsion-free $H^\#$ iff $G$ contains no copies of
$V_\omega$. One can consider also ``local'' forms, e.g., about continuous
maps $B_\omega^\#\to H^\#$ (or $V_\omega^\#\to H^\#$ for a particular {\em
fixed} $p$).

To introduce the next stronger form of the Straightening Law we isolate
first the following corollary of Theorem \ref{lastT}. 

\begin{cor} 
If $G$ is a torsion group of size $>\beth_{\omega}$, then every
continuous map $G^\#\to H^\#$ is linear on some infinite set $A\subseteq 
G$. 
\end{cor}

Clearly, ``torsion'' can be replaced by $|\mbox{tor}(G)|>\beth_\omega$. 

\begin{question} 
Can ``torsion'' be completely removed in the above corollary? 
\end{question}

In other words, does the corollary work with groups with ``small'' torsion
part? Therefore, the second form (we formulate it as a question) is
determined by {\em torsion-free domain}.

\begin{question} 
If $G$ is a sufficiently large torsion-free group, is then every
continuous map $G^\#\to H^\#$ linear on some infinite set $A\subseteq G$? 
\end{question}

Clearly, ``torsion-free'' can be replaced by ``free abelian'', i.e., {\em
is every continuous map $(\bigoplus_\kappa\Z)^\#\to H^\#$ linear on some
infinite set $A\subseteq \bigoplus_\kappa\Z$} ? In particular, the
following question is left open.

\begin{question} 
Does there exist a continuous injective map $\Z^\#\to T^\#$, with $T$
torsion abelian group? 
\end{question}

One can always take $T=\bigoplus_\omega \Q/\Z$. 

\subsubsection*{}

Now we shall discuss the third stronger form of the Straightening Law --
when restriction to {\em subgroups} of continuous maps (in the Bohr
topology) are linear. According to Theorem \ref{ST}, if $\kappa>\beth_3$
and $f:B_\kappa^\#\to H^\#$ is a continuous map, then $f$ is linear on
infinitely many non-trivial subgroups of $B_\kappa$. Namely, if $f$ is
linear on $[S]^4$ for an infinite $S\subseteq \kappa$, then $f$ is linear
on $A=S+S+S+S=\{0\}\cup [S]^2\cup [S]^4=\overline{[S]^4}$, so $f$ is
linear on every cyclic subgroup $span(a)$, $a\in A$, as well as on any
2-generated subgroup of $B_\kappa$ with generators from $S$. It is clear,
that taking a larger $\kappa$ (depending on $k$), one can obtain linearity
on $A=\{0\}\cup [S]^2\cup \ldots \cup [S]^{2k}=\overline{[S]^{2k}}$ for
any fixed $k\in \N$. Then $f$ is linear on every $2k$-generated subgroup
of $B_\kappa$ with generators from $S$. This leaves open the following:

\begin{question} 
Is $f$ linear on some {\em infinite subgroup} of $B_\kappa$ if $\kappa$
is sufficiently large? 
\end{question}

\begin{question} 
When every continuous map $G^\#\to H^\#$ is linear on some infinite {\em
subgroup} of $G$? 
\end{question}

Possible restraints to impose on $G$ are: ``large'', Boolean, 
$G\cong V_\kappa$ for some large $\kappa$. 

\begin{question} 
Do there exist continuous maps $f: G^\#\to H^\#$ with uncountable $G$,
that are not linear on any infinite {\em subgroup} of $G$? 
\end{question} 

Clearly, one has to rule out maps as $j: (\Q/\Z)^\#\hookrightarrow \Q^\#$
by asking the existence of infinite subgroups $A$ of $G$ with 
$Hom(A,H) \ne 0$.

\section{Retracts of $G^\#$} 

van Douwen posed also the following questions about retracts of $G^\#$:

\begin{question}\label{retr_subgr} {\rm \cite[Question 4.12]{vD}} 
Is every subgroup $H$ of a group $G^{\#}$ a retract of $G^{\#}$? 
\end{question} 

Let us recall that every subgroup of $G^\#$ is closed. 

\begin{question}\label{subset} {\rm \cite[Question 4.13]{vD}} 
Is every countable closed subset of $G^{\#}$ a retract of $G^{\#}$?
\end{question} 

\begin{remark} 
If a closed {\em relatively discrete} set $R\subseteq G^\#$ is a retract
of $G^{\#}$, then $R$ is {\em countable}. 
Indeed, if $r:G^\# \to R$ is a retract, then $r^{-1}(R)$ is disjoint union
of $|R|$ open pairwise disjoint sets of $G^\#$. 
Since $c(G)=\omega$, this implies $R$ is countable. 
It remains to note that if $G$ is uncountable, then $G^\#$ admits a
relatively discrete subset $R$ of size $|G|$ (by Theorem 
\ref{vD_Thm}). 
This explains why {\em countable} is imposed in Question
\ref{retr_subgr}. 
\end{remark}

A negative answer to the second question was obtained by Gladdines
\cite{gladd} in 1995. Now we offer a new much shorter proof of her
theorem based on the non-homeomorphism theorems and the better knowledge
of the continuous maps in the Bohr topology obtained after 1995.

\begin{theorem} {\rm\cite{gladd}} 
$\mathcal{D}_\kappa=\{0\}\cup [\omega]^2$ is not a retract of
$B_\omega^\#$.
\end{theorem}

\begin{proof} 
Assume that $r:B_\omega^\#\to \{0\}\cup [\omega]^2$ is a
retraction. Let $H= \bigoplus_\omega \Z_3$, let $\langle
e_n:n<\omega\rangle$ be the canonical base of $H$ and let $\mu:\{0\}\cup
[\omega]^2\to H^\#$ be the continuous map defined by $\mu(0)=0$ and
$\mu(n,m)=e_n-e_m$ for $n<m<\omega$. Then $\pi=\mu\circ r:B_\omega^\#\to
H^\#$ is continuous with $|\supp \pi(x)|\leq 2$ for all $x\in
B_\omega$, so there exists an infinite subset $S\subseteq \omega$ such
that $\pi$ vanishes on $[S]^4$ (cf.\ Theorem \ref{K_THM} or Remark
\ref{REM} (a)), hence also on $[S]^2\subseteq \overline {[S]^4}$. A
contradiction, since $\pi\upharpoonright [S]^2=\mu \upharpoonright [S]^2$
is injective (being $r\upharpoonright [S]^2=id_{[S]^2}$).
\end{proof} 

\subsection{Subgroups as retracts}\label{SecCHT}

Question \ref{subset} is still open. Here we offer a comment on what was
done so far.

In $B_\kappa$ every subgroup $L$ splits off $B_\kappa=L\times K$, hence
$L$ is a topological group retract in $B_\kappa^\#=L^\#\times K^\#$. The
same applies to $V_\kappa^\#$. In general, if $H$ has finite index in
$G$, then $H$ is clopen in $G^{\#}$, hence $H$ is a retract of $G^{\#}$
(for other instances see \cite{T1}). 

A substantial contribution towards a solution to this problem was given by
Comfort, Hern\'andez and Trigos-Arrieta \cite{CHT3} by the introduction of
an important notion that helps to understand better the nature of the
retract problem in the case of subgroups. They call a subgroup $H$ of an
abelian group $G$ a {\em ccs-subgroup} if the the natural map $\varphi:
G^\#\to G^\#/H$ has a cross section $\Gamma: G/H \to G$ that is {\em
continuous } in the Bohr topology of $G$ and $G/H$. It is easy to see that
in such a case $G^\#$ is homeomorphic to $(G/H)^\#\times H^\#$. In
particular, $H$ is a retract of $G^\#$ and $(G/H)^\#$ embeds into $G^\#$.
The existence of a Bohr-continuous cross section $\Gamma: G/H \to G$ is
equivalent to the existence of a retraction $G^\#\to H^\#$ such that
$r(x+h)=r(x)+h$ for every $h\in H$ and every $x\in G$ (so that $r$ is a
{\em linear} retraction \cite{CHT3}). In this sense the study of
ccs-subgroups, proposed in \cite{CHT3}, is a very natural modification of
van Douwen's Question \ref{retr_subgr}.

Comfort, Hern\'andez and Trigos-Arrieta \cite{CHT3} introduced the
class ${\bf ACCS}(\#)$ of abelian groups $H$ that are ccs-subgroups of any
abelian group that contains them and they showed that $H\in {\bf
ACCS}(\#)$ iff $H$ is a ccs-subgroup of its the divisible hull. They
showed that the class ${\bf ACCS}(\#)$ is closed under finite products and
contains all finitely generated groups (and, of course, all divisible
groups).

Here are some further examples of groups from ${\bf ACCS}(\#)$:

\begin{example}
\mbox{}
\begin{itemize}
\item[(a)] 
{\rm (\cite{CHT3})} $\Z^n\in {\bf ACCS}(\#)$ and consequently, 
$(\Q^n/\Z^n)^\#$ embeds into $(\Q^n)^\#$ for every $n<\omega$. We repeat
this example in view of its importance.
\item[(b)] 
{\rm (\cite{CHT3})} For every prime $p$ the group 
$\bigoplus_\omega\Z_{p}$ is not a ccs-subgroup of
$\bigoplus_\omega\Z_{p^2}$ (so the group $\bigoplus_\omega\Z_{p}$
does not belong to ${\bf ACCS}(\#)$).
\item[(c)] 
{\rm (\cite{DHH})} Let $\pi$ be an arbitrary non-empty set of prime
numbers. 
\begin{itemize}
\item[(i)] 
The additive group $H_\pi$ of the subring of $\Q$ generated by the set
$\{1/p:p\in \pi\}$ belongs to ${\bf ACCS}(\#)$ (for a singleton
$\pi=\{p\}$ this was mentioned also in \cite{CHT3}). By varying $\pi$
this gives $\co$ many pairwise non-isomorphic rank-one torsion-free
(reduced) groups in ${\bf ACCS}(\#)$. 
\item[(ii)] 
The product $\prod_{p\in \pi}\J_p$ belongs to ${\bf ACCS}(\#)$, where
$\J_p$ denotes the group of $p$-adic integers. By varying $\pi$ this
gives $\co$ many pairwise non-isomorphic reduced groups of size $\co$.
\end{itemize}
\end{itemize}
\end{example} 

Some new restraints for the groups from ${\bf ACCS}(\#)$ are given in 
\cite{DHH}. 
They give an upper bound for the size of the reduced groups in ${\bf
ACCS}(\#)$ (so that the reduced groups in ${\bf ACCS}(\#)$ form a
set) and entails that large powers may belong to ${\bf ACCS}(\#)$ only
if they are divisible. Under the assumption that the Strong Straightening
Law holds, every reduced group $H\in {\bf ACCS}(\#)$ is almost
torsion-free (so its torsion part is countable) and $|H|\leq \co$
\cite{DHH}. 

\begin{question}
\mbox{} 
\begin{itemize}
\item[(a)] 
{\rm (\cite{DHH})} Is the subgroup $H$ of $\Q$ spanned by all fractions
$1/p$, with $p$ prime, a retract (ccs-subgroup) of $\Q^\#$ ?
\item[(b)] 
{\rm (\cite{CHT3})} Is the subgroup $\bigoplus_\omega\Z_p$ of 
$G=\bigoplus_\omega\Z_{p^2}$ a retract of $G^\#$?
\end{itemize}
\end{question}

``Yes'' to item (a) is equivalent to $\bigoplus_p \Z_p \in {\bf ACCS}(\#)$
\cite{DHH}. 
Item (b) is open even for $p=2$. 

\providecommand{\bysame}{\leavevmode\hbox to3em{\hrulefill}\thinspace}
\providecommand{\MR}{\relax\ifhmode\unskip\space\fi MR }
\providecommand{\MRhref}[2]{%
  \href{http://www.ams.org/mathscinet-getitem?mr=#1}{#2}
}
\providecommand{\href}[2]{#2}

\end{document}